# Biring of Matrices

Aleks Kleyn

Abstract. Matrices allow two products linked by transpose. Biring is algebra which defines on the set two correlated structures of the ring. According to each product we can extend the definition of a quasideterminant given in [1, 2] and introduce two different types of a quasideterminant.

Contents

## 1. Concept of Generalized Index

Studying tensor calculus we start from studying univalent covariant and contravariant tensors. In spite on difference of properties both these objects are elements of respective vector spaces. Suppose we introduce a generalized index according to the rule $a^i = a^i$, $b^i = b^{\cdot}{}_{i}^{-}$. Then we see that these tensors have the similar behavior. For instance, the transformation of a covariant tensor gets form

$$b'^i = b'^{\cdot}{}_{i}^{-} = f^{\cdot}{}_{i}^{-}{}^{\cdot}{}_{-}^{j} b^{\cdot}{}_{j}^{-} = f^i_j b^j$$

This similarity goes as far as we need because tensors also form vector space.

These observations of the similarity between properties of covariant and contravariant tensors lead us to the concept of generalized index.[1] We will use the symbol $\cdot$ in front of a generalized index when we need to describe its structure. I put the sign $'-'$ in place of the index whose position was changed. For instance, if an original term was $a_{ij}$ I will use notation $a_i{}^j_-$ instead of notation $a^j_i$.

Aleks_Kleyn@MailAPS.org.
http://sites.google.com/site/AleksKleyn/.
http://arxiv.org/a/kleyn_a_1.
http://AleksKleyn.blogspot.com/.

[1] Set of commands to use generalized index as well other commands that I use in this and following papers you can see in the file Commands.tex.

Even though the structure of a generalized index is arbitrary we assume that there exists a one-to-one map of the interval of positive integers 1, ..., $n$ to the range of index. Let $i$ be the range of the index $i$. We denote the power of this set by symbol $|i|$ and assume that $|i| = n$. If we want to enumerate elements $a_i$ we use notation $a_1$, ..., $a_n$.

Representation of coordinates of a vector as a matrix allows making a notation more compact. The question of the presentation of vector as a row or a column of the matrix is just a question of convention. We extend the concept of generalized index to elements of the matrix. A matrix is a two dimensional table, the rows and columns of which are enumerated by generalized indexes. To represent a matrix we will use one of the following forms:

**Standard representation:** in this case we write elements of matrix $A$ as $A^a_b$.

**Alternative representation:** in this case we write elements of matrix $A$ as ${}^aA_b$ or ${}_bA^a$.

Since we use generalized index, we cannot tell whether index $a$ of matrix enumerates rows or columns until we know the structure of index.

We could use notation $*$-column and $*$-row which is more close to our custom. However as we can see bellow the form of presentation of matrix is not important for us. To make sure that notation offered below is consistent with the traditional we will assume that the matrix is presented in the form

$$A = \begin{pmatrix} {}^1A_1 & ... & {}^1A_n \\ ... & ... & ... \\ {}^mA_1 & ... & mA_n \end{pmatrix}$$

**Definition 1.1.** I use the following names and notation for different minors of the matrix $A$

$A_a$ **:** $^*$**-row** with the index $a$ is generalization of a column of a matrix. The upper index enumerates elements of $^*$-rows and the lower index enumerates $^*$-rows.

$A_T$ **:** the minor obtained from $A$ by selecting $^*$-rows with an index from the set $T$

$A_{[a]}$ **:** the minor obtained from $A$ by deleting $^*$-row $A_a$

$A_{[T]}$ **:** the minor obtained from $A$ by deleting $^*$-rows with an index from the set $T$

${}^bA$ **:** ${}_*$**-row** with the index $b$ is generalization of a row of a matrix. The lower index enumerates elements of ${}_*$-rows and the upper index enumerates ${}_*$-rows.

${}^SA$ **:** the minor obtained from $A$ by selecting ${}_*$-rows with an index from the set $S$

${}^{[b]}A$ **:** the minor obtained from $A$ by deleting ${}_*$-row ${}^bA$

${}^{[S]}A$ **:** the minor obtained from $A$ by deleting ${}_*$-rows with an index from the set $S$

□

*Remark* 1.2. We will combine the notation of indexes. Thus ${}^bA_a$ is $1 \times 1$ minor. The same time this is the notation for a matrix element. This allows an identifying of $1 \times 1$ matrix and its element. The index $a$ is number of $^*$-row of matrix and the index $b$ is number of ${}_*$-rows of matrix. □

Each form of the notation of a matrix has its own advantages. The standard notation is more natural when we study matrix theory. The alternative form of the notation makes expressions in the theory of vector spaces more clear. Extending the alternative notation of indexes to arbitrary tensors we can better understand an interaction of different geometric objects. Using the duality principle (theorem 2.14) improves our expressivity.

*Remark* 1.3. We can read symbol ${}^*$- as $c$- and symbol ${}_*$- as $r$- creating this way names $c$-**row** and $r$-**row**. Further we extend this rule to other objects of linear algebra. I will use this convention designing index. □

Since transpose of the matrix exchanges ${}_*$-rows and ${}^*$-rows we get equation

$$ {}_i(A^T)^j = {}^iA_j \tag{1.1} $$

*Remark* 1.4. As we can see from the equation (1.1), it is not important for us the choice of a side to place a number of ${}_*$-row and and the choice of a side to place a number of ${}^*$-row. This is due to the fact that we can enumerate the elements of matrix in different ways. If we want to show the numbers of ${}_*$-row and ${}^*$-row according to the definition 1.1, then the equation (1.1) has form

$$ {}^j(A^T)_i = {}^iA_j $$

In standard representation, the equation (1.1) has form

$$ (A^T)^j_i = A^i_j $$

□

We call matrix[2]

$$ \mathcal{H}A = ({}^j\mathcal{H}A_i) = ((^{\cdot}{}_{i}^{-}A^{\cdot}{}_{-}^{j})^{-1}) \tag{1.2} $$

**Hadamard inverse of matrix** $A = ({}_bA^a)$ ([2]-page 4).

I will use the Einstein convention about sums. This means that when an index is present in an expression twice and a set of index is known, I have the sum over this index. If needed to clearly show set of index, I will do it. Also, in this paper I will use the same root letter for a matrix and its elements.

We will study matrices elements of which belong to division ring $D$. We will also keep in mind that instead of division ring $D$ we may write in text field $F$. We will clearly write field $F$ in case when commutativity creates new details. We will denote by 1 identity element of division ring $D$.

Let $I$, $|I| = n$ be a set of indexes. We introduce the **Kronecker symbol**

$$ \delta^i_j = \begin{cases} 1 & i = j \\ 0 & i \ne j \end{cases} \qquad i, j \in I \tag{1.3} $$

---

[2] The notation $(^{\cdot}{}_{i}^{-}A^{\cdot}{}_{-}^{j})^{-1}$ means that we exchange rows and columns in Hadamard inverse. We can formally write this expression in following form

$$ (^{\cdot}{}_{i}^{-}A^{\cdot}{}_{-}^{j})^{-1} = \frac{1}{{}^iA_j} $$

## 2. Biring

We consider matrices whose elements belong to division ring $D$.

The product of matrices is associated with the product of homomorphisms of vector spaces over field. According to the custom the product of matrices $A$ and $B$ is defined as product of ${}_*$-rows of the matrix $A$ and ${}^*$-rows of the matrix $B$. Conventional character of this definition becomes evident when we put attention that ${}_*$-row of the matrix $A$ may be a column of this matrix. In such case we multiply columns of the matrix $A$ over rows of the matrix $B$. Thus we can define two products of matrices. To distinguish between these products we introduced a new notation.[3]

**Definition 2.1.** **${}_*{}^*$-product of matrices** $A$ and $B$ has form

$$\left\{\begin{array}{rcl}A_*{}^*B&=&({}^aA_c\ {}^cB_b)\\ {}^a(A_*{}^*B)_b&=&{}^aA_c\ {}^cB_b\end{array}\right. \tag{2.1}$$

and can be expressed as product of a ${}_*$-row of matrix $A$ over a ${}^*$-row of matrix $B$.[4] □

**Definition 2.2.** **${}^*{}_*$-product of matrices** $A$ and $B$ has form

$$\left\{\begin{array}{rcl}A^*{}_*B&=&({}_aA^c\ {}_cB^b)\\ {}_a(A^*{}_*B)^b&=&{}_aA^c\ {}_cB^b\end{array}\right. \tag{2.2}$$

and can be expressed as product of a ${}^*$-row of matrix $A$ over a ${}_*$-row of matrix $B$.[5] □

*Remark* 2.3. We will use symbol ${}_*{}^*$- or ${}^*{}_*$- in name of properties of each product and in the notation. According to remark 1.3 we can read symbols ${}_*{}^*$ and ${}^*{}_*$ as *rc* -product and *cr*-product. This rule we extend to following terminology. □

*Remark* 2.4. Just as in remark 1.4, I want to draw attention to the fact that I change the numbering of elements of the matrix. If we want to show the numbers of ${}_*$-row and ${}^*$-row according to the definition 1.1, then the equation (2.2) has form

$$ {}^b(A^*{}_*B)_a={}^cA_a\ {}^bB_c \tag{2.3}$$

However the format of the equation (2.3) is unusual. □

---

[3] In order to keep this notation consistent with the existing one we assume that we have in mind ${}_*{}^*$-product when no clear notation is present.

[4] In alternative form operation consists from two symbols $*$ which we put in the place of index which participate in sum. In standard notation we write operation as

$$\left\{\begin{array}{rcl}A_*{}^*B&=&(A^a_cB^c_b)\\ (A_*{}^*B)^a_b&=&A^a_cB^c_b\end{array}\right.$$

and can be construed as symbolic notation

$$A_*{}^*B=A_*B^*$$

where we write symbol $*$ on place of index which participate in sum.

[5] In alternative form operation consists from two symbols $*$ which we put in the place of index which participate in sum. In standard notation we write operation as

$$\left\{\begin{array}{rcl}A^*{}_*B&=&(A^c_a\ B^b_c)\\ (A^*{}_*B)^a_b&=&A^c_a\ B^b_c\end{array}\right.$$

and can be construed as symbolic notation

$$A^*{}_*B=A^*B_*$$

where we write symbol $*$ on place of index which participate in sum.

Set of $n \times n$ matrices is closed relative ${}_*{}^*$-product and ${}^*{}_*$-product as well relative sum which is defined by rule

$$(A+B)_a^b = A_a^b + B_a^b$$

**Theorem 2.5.**

$$(A_*{}^*B)^T = A^T{}^*{}_*B^T \tag{2.4}$$

*Proof.* The chain of equations

$$\begin{aligned} {}_a((A_*{}^*B)^T)^b &= {}^a(A_*{}^*B)_b \\ &= {}^aA_c\ {}^cB_b \\ &= {}_a(A^T)^c{}_c(B^T)^b \\ &= {}_a((A^T)^*{}_*(B^T))^b \end{aligned} \tag{2.5}$$

follows from (1.1), (2.1) and (2.2). The equation (2.4) follows from (2.5). □

Matrix $\delta = (\delta_a^c)$ is identity for both products.

**Definition 2.6.** $\mathcal{A}$ is a **biring** if we defined on $\mathcal{A}$ an unary operation, say transpose, and three binary operations, say ${}_*{}^*$-product, ${}^*{}_*$-product and sum, such that

- ${}_*{}^*$-product and sum define structure of ring on $\mathcal{A}$
- ${}^*{}_*$-product and sum define structure of ring on $\mathcal{A}$
- both products have common identity $\delta$
- products satisfy equation

$$(A_*{}^*B)^T = A^T{}^*{}_*B^T \tag{2.6}$$

- transpose of identity is identity

$$\delta^T = \delta \tag{2.7}$$

- double transpose is original element

$$(A^T)^T = A \tag{2.8}$$

□

**Theorem 2.7.**

$$(A^*{}_*B)^T = (A^T)_*{}^*(B^T) \tag{2.9}$$

*Proof.* We can prove (2.9) in case of matrices the same way as we proved (2.6). However it is more important for us to show that (2.9) follows directly from (2.6).

Applying (2.8) to each term in left side of (2.9) we get

$$(A^*{}_*B)^T = ((A^T)^T{}^*{}_*(B^T)^T)^T \tag{2.10}$$

From (2.10) and (2.6) it follows that

$$(A^*{}_*B)^T = ((A^T{}_*{}^*B^T)^T)^T \tag{2.11}$$

(2.9) follows from (2.11) and (2.8). □

**Definition 2.8.** We introduce ${}_*{}^*$**-power** of element $A$ of biring $\mathcal{A}$ using recursive definition

$$A^{0_*{}^*} = \delta \tag{2.12}$$

$$A^{n_*{}^*} = A^{n-1_*{}^*}{}_*{}^*A \tag{2.13}$$

□

**Definition 2.9.** We introduce **${}^*{}_*$-power** of element $A$ of biring $\mathcal{A}$ using recursive definition

$$A^{0^*{}_*}=\delta \tag{2.14}$$

$$A^{n^*{}_*}=A^{n-1^*{}_*}{}^*{}_*A \tag{2.15}$$

□

**Theorem 2.10.**

$$(A^T)^{n_*{}^*}=(A^{n^*{}_*})^T \tag{2.16}$$

$$(A^T)^{n^*{}_*}=(A^{n_*{}^*})^T \tag{2.17}$$

*Proof.* We proceed by induction on $n$.

For $n=0$ the statement immediately follows from equations (2.12), (2.14), and (2.7).

Suppose the statement of theorem holds when $n=k-1$

$$(A^T)^{n-1_*{}^*}=(A^{n-1^*{}_*})^T \tag{2.18}$$

It follows from (2.13) that

$$(A^T)^{k_*{}^*}=(A^T)^{k-1_*{}^*}{}_*{}^*A^T \tag{2.19}$$

It follows from (2.19) and (2.18) that

$$(A^T)^{k_*{}^*}=(A^{k-1^*{}_*})^T{}_*{}^*A^T \tag{2.20}$$

It follows from (2.20) and (2.9) that

$$(A^T)^{k_*{}^*}=(A^{k-1^*{}_*}{}^*{}_*A)^T \tag{2.21}$$

(2.16) follows from (2.19) and (2.15).

We can prove (2.17) by similar way. □

**Definition 2.11.** Element $A^{-1_*{}^*}$ of biring $\mathcal{A}$ is **${}_*{}^*$-inverse element** of element $A$ if

$$A{}_*{}^*A^{-1_*{}^*}=\delta \tag{2.22}$$

Element $A^{-1^*{}_*}$ of biring $\mathcal{A}$ is **${}^*{}_*$-inverse element** of element $A$ if

$$A{}^*{}_*A^{-1^*{}_*}=\delta \tag{2.23}$$

□

**Theorem 2.12.** *Suppose element $A\in\mathcal{A}$ has ${}_*{}^*$-inverse element. Then transpose element $A^T$ has ${}^*{}_*$-inverse element and these elements satisfy equation*

$$(A^T)^{-1^*{}_*}=(A^{-1_*{}^*})^T \tag{2.24}$$

*Suppose element $A\in\mathcal{A}$ has ${}^*{}_*$-inverse element. Then transpose element $A^T$ has ${}_*{}^*$-inverse element and these elements satisfy equation*

$$(A^T)^{-1_*{}^*}=(A^{-1^*{}_*})^T \tag{2.25}$$

*Proof.* If we get transpose of both side (2.22) and apply (2.7) we get

$$(A_*{}^*A^{-1_*{}^*})^T = \delta^T = \delta$$

Applying (2.6) we get

$$\delta = A^T{}^*{}_*(A^{-1_*{}^*})^T \tag{2.26}$$

(2.24) follows from comparison (2.23) and (2.26).

We can prove (2.25) similar way. □

Theorems 2.5, 2.7, 2.10, and 2.12 show that some kind of duality exists between ${}_*{}^*$-product and ${}^*{}_*$-product. We can combine these statements.

**Theorem 2.13 (duality principle for biring).** *Let $\mathfrak{A}$ be true statement about biring $\mathcal{A}$. If we exchange the same time*

- *$A \in \mathcal{A}$ and $A^T$*
- *${}_*{}^*$-product and ${}^*{}_*$-product*

*then we soon get true statement.*

**Theorem 2.14 (duality principle for biring of matrices).** *Let $\mathcal{A}$ be biring of matrices. Let $\mathfrak{A}$ be true statement about matrices. If we exchange the same time*

- *${}^*$-rows and ${}_*$-rows of all matrices*
- *${}_*{}^*$-product and ${}^*{}_*$-product*

*then we soon get true statement.*

*Proof.* This is the immediate consequence of the theorem 2.13. □

*Remark* 2.15. We execute operations in expression

$$A_*{}^*B_*{}^*C$$

from left to right. However we can execute product from right to left. In custom notation this expression is

$$C^*{}_*B^*{}_*A$$

We follow the rule that to write power from right of expression. If we use standard representation, then we write indexes from right of expression. If we use alternative representation, then we read indexes in the same order as symbols of operation and root letters. For instance, let original expression be like

$$A^{-1_*{}^*}{}_*{}^*B_a$$

Then expression which we reed from right to left is like

$$B_a{}^*{}_*A^{-1^*{}_*}$$

in standard representation or

$${}_aB^*{}_*A^{-1^*{}_*}$$

in alternative representation.

Suppose we established the order in which we write indexes. Then we state that we read an expression from top to bottom reading first upper indexes, then lower ones. We assume that this is standard form of reading. We can read this expression from bottom to top. We extend this rule stating that we read symbols of operation in the same order as indexes. For instance, if we read expression

$$A^a{}_*{}^*B^{-1_*{}^*} = C^a$$

from bottom to top, then we can write this expression in standard form

$$A_a{}^*{}_*B^{-1^*{}_*} = C_a$$

According to the duality principle if we can prove one statement then we can prove other as well. □

**Theorem 2.16.** *Let matrix $A$ have ${}_*{}^*$-inverse matrix. Then for any matrices $B$ and $C$ equation*

$$B = C \tag{2.27}$$

*follows from the equation*

$$B_*{}^*A = C_*{}^*A \tag{2.28}$$

*Proof.* Equation (2.27) follows from the equation (2.28) if we multiply both parts of the equation (2.28) over $A^{-1_*{}^*}$. □

## 3. Quasideterminant

**Theorem 3.1.** *Suppose $n\times n$ matrix $A$ has ${}_*{}^*$-inverse matrix.*[6] *Then $k\times k$ minor of ${}_*{}^*$-inverse matrix satisfy*

$$\left({}^I(A^{-1_*{}^*})_J\right)^{-1_*{}^*} = {}^JA_I - {}^JA_{[I]}{}_*{}^*\left({}^{[J]}A_{[I]}\right)^{-1_*{}^*}{}_*{}^*{}^{[J]}A_I \tag{3.1}$$

*Proof.* Definition (2.22) of ${}_*{}^*$-inverse matrix leads to system of linear equations

$${}^{[J]}A_{[I]}{}_*{}^*{}^{[I]}(A^{-1_*{}^*})_J + {}^{[J]}A_I{}_*{}^*{}^I(A^{-1_*{}^*})_J = 0 \tag{3.2}$$

$${}^JA_{[I]}{}_*{}^*{}^{[I]}(A^{-1_*{}^*})_J + {}^JA_I{}_*{}^*{}^I(A^{-1_*{}^*})_J = \delta \tag{3.3}$$

We multiply (3.2) by $\left({}^{[J]}A_{[I]}\right)^{-1_*{}^*}$

$${}^{[I]}(A^{-1_*{}^*})_J + \left({}^{[J]}A_{[I]}\right)^{-1_*{}^*}{}_*{}^*{}^{[J]}A_I{}_*{}^*{}^I(A^{-1_*{}^*})_J = 0 \tag{3.4}$$

Now we can substitute (3.4) into (3.3)

$$-{}^JA_{[I]}{}_*{}^*\left({}^{[J]}A_{[I]}\right)^{-1_*{}^*}{}_*{}^*{}^{[J]}A_I{}_*{}^*{}^I(A^{-1_*{}^*})_J + {}^JA_I{}_*{}^*{}^I(A^{-1_*{}^*})_J = \delta \tag{3.5}$$

(3.1) follows from (3.5). □

**Corollary 3.2.** *Suppose $n\times n$ matrix $A$ has ${}_*{}^*$-inverse matrix. Then elements of ${}_*{}^*$-inverse matrix satisfy to the equation*[2]

$${}^i(A^{-1_*{}^*})_j = \left({}^jA_i - {}^jA_{[i]}{}_*{}^*\left({}^{[j]}A_{[i]}\right)^{-1_*{}^*}{}_*{}^*{}^{[j]}A_i\right)^{-1} \tag{3.6}$$

$${}^j\left(\mathcal{H}A^{-1_*{}^*}\right)_i = {}^jA_i - {}^jA_{[i]}{}_*{}^*\left({}^{[j]}A_{[i]}\right)^{-1_*{}^*}{}_*{}^*{}^{[j]}A_i \tag{3.7}$$

□

[6] This statement and its proof are based on statement 1.2.1 from [1] (page 8) for matrix over free division ring.

**Example 3.3.** Consider matrix

$$\begin{pmatrix} {}^1A_1 & {}^1A_2 \\ {}^2A_1 & {}^2A_2 \end{pmatrix}$$

According to (3.6)

$$ {}^1(A^{-1_*{}^*})_1 = ({}^1A_1 - {}^1A_2({}^2A_2)^{-1}\,{}^2A_1)^{-1} \tag{3.8}$$

$$ {}^2(A^{-1_*{}^*})_1 = ({}^1A_2 - {}^1A_1({}^2A_1)^{-1}\,{}^2A_2)^{-1} \tag{3.9}$$

$$ {}^1(A^{-1_*{}^*})_2 = ({}^2A_1 - {}^2A_2({}^1A_2)^{-1}\,{}^1A_1)^{-1} \tag{3.10}$$

$$ {}^2(A^{-1_*{}^*})_2 = ({}^2A_2 - {}^2A_1({}^1A_1)^{-1}\,{}^1A_2)^{-1} \tag{3.11}$$

Consider the product of matrices

$$C = \begin{pmatrix} ({}^1A_1 - {}^1A_2({}^2A_2)^{-1}\,{}^2A_1)^{-1} & ({}^2A_1 - {}^2A_2({}^1A_2)^{-1}\,{}^1A_1)^{-1} \\ ({}^1A_2 - {}^1A_1({}^2A_1)^{-1}\,{}^2A_2)^{-1} & ({}^2A_2 - {}^2A_1({}^1A_1)^{-1}\,{}^1A_2)^{-1} \end{pmatrix} {}_*{}^* \begin{pmatrix} {}^1A_1 & {}^1A_2 \\ {}^2A_1 & {}^2A_2 \end{pmatrix}$$

From direct calculations, it follows that

$$\begin{aligned}
{}^1C_1 &= ({}^1A_1 - {}^1A_2({}^2A_2)^{-1}\,{}^2A_1)^{-1}\,{}^1A_1 + ({}^2A_1 - {}^2A_2({}^1A_2)^{-1}\,{}^1A_1)^{-1}\,{}^2A_1 \\
&= ({}^1A_2(({}^1A_2)^{-1}\,{}^1A_1 - ({}^2A_2)^{-1}\,{}^2A_1))^{-1}\,{}^1A_1 \\
&+ ({}^2A_2(({}^2A_2)^{-1}\,{}^2A_1 - ({}^1A_2)^{-1}\,{}^1A_1))^{-1}\,{}^2A_1 \\
&= (({}^1A_2)^{-1}\,{}^1A_1 - ({}^2A_2)^{-1}\,{}^2A_1)^{-1}\ ({}^1A_2)^{-1}\,{}^1A_1 \\
&+ (({}^2A_2)^{-1}\,{}^2A_1 - ({}^1A_2)^{-1}\,{}^1A_1)^{-1}\ ({}^2A_2)^{-1}\,{}^2A_1 \\
&= 1
\end{aligned}$$

$$\begin{aligned}
{}^1C_2 &= ({}^1A_1 - {}^1A_2({}^2A_2)^{-1}\,{}^2A_1)^{-1}\,{}^1A_2 + ({}^2A_1 - {}^2A_2({}^1A_2)^{-1}\,{}^1A_1)^{-1}\,{}^2A_2 \\
&= ({}^1A_2(({}^1A_2)^{-1}\,{}^1A_1 - ({}^2A_2)^{-1}\,{}^2A_1))^{-1}\,{}^1A_2 \\
&+ ({}^2A_2(({}^2A_2)^{-1}\,{}^2A_1 - ({}^1A_2)^{-1}\,{}^1A_1))^{-1}\,{}^2A_2 \\
&= (({}^1A_2)^{-1}\,{}^1A_1 - ({}^2A_2)^{-1}\,{}^2A_1)^{-1}\ ({}^1A_2)^{-1}\,{}^1A_2 \\
&+ (({}^2A_2)^{-1}\,{}^2A_1 - ({}^1A_2)^{-1}\,{}^1A_1)^{-1}\ ({}^2A_2)^{-1}\,{}^2A_2 \\
&= 0
\end{aligned}$$

$$\begin{aligned}
{}^2C_1 &= ({}^1A_2 - {}^1A_1({}^2A_1)^{-1}\,{}^2A_2)^{-1}\,{}^1A_1 + ({}^2A_2 - {}^2A_1({}^1A_1)^{-1}\,{}^1A_2)^{-1}\,{}^2A_1 \\
&= ({}^1A_1(({}^1A_1)^{-1}\,{}^1A_2 - ({}^2A_1)^{-1}\,{}^2A_2)^{-1}\,{}^1A_1 \\
&+ ({}^2A_1(({}^2A_1)^{-1}\,{}^2A_2 - ({}^1A_1)^{-1}\,{}^1A_2))^{-1}\,{}^2A_1 \\
&= (({}^1A_1)^{-1}\,{}^1A_2 - ({}^2A_1)^{-1}\,{}^2A_2)^{-1}\ ({}^1A_1)^{-1}\,{}^1A_1 \\
&+ (({}^2A_1)^{-1}\,{}^2A_2 - ({}^1A_1)^{-1}\,{}^1A_2)^{-1}\ ({}^2A_1)^{-1}\,{}^2A_1 \\
&= 0
\end{aligned}$$

$$
\begin{aligned}
{}^2C_2 &= ({}^1A_2 - {}^1A_1({}^2A_1)^{-1}\ {}^2A_2)^{-1}\ {}^1A_2 + ({}^2A_2 - {}^2A_1({}^1A_1)^{-1}\ {}^1A_2)^{-1}\ {}^2A_2\\
&= ({}^1A_1(({}^1A_1)^{-1}\ {}^1A_2 - ({}^2A_1)^{-1}\ {}^2A_2)^{-1}\ {}^1A_2\\
&+ ({}^2A_1(({}^2A_1)^{-1}\ {}^2A_2 - ({}^1A_1)^{-1}\ {}^1A_2))^{-1}\ {}^2A_2\\
&= ((({}^1A_1)^{-1}\ {}^1A_2 - ({}^2A_1)^{-1}\ {}^2A_2)^{-1}\ ({}^1A_1)^{-1}\ {}^1A_2\\
&+ (({}^2A_1)^{-1}\ {}^2A_2 - ({}^1A_1)^{-1}\ {}^1A_2)^{-1}\ ({}^2A_1)^{-1}\ {}^2A_2\\
&= 1
\end{aligned}
$$

□

According to [1], page 3 we do not have an appropriate definition of a determinant for a division ring. However, we can define a quasideterminant which finally gives a similar picture. In definition below we follow definition [1]-1.2.2.

**Definition 3.4.** $\mathbf{\binom{j}{i}}$**-$_*{}^*$-quasideterminant** of $n\times n$ matrix $A$ is formal expression[2]

$$
{}^j\det(A,{}_*{}^*)_i = {}^j\left(\mathcal{H}A^{-1_*{}^*}\right)_i \tag{3.12}
$$

According to the remark 1.2 we can get $\binom{j}{i}$-$_*{}^*$-quasideterminant as an element of the matrix $\det(a,{}_*{}^*)$ which we call $_*{}^*$**-quasideterminant**. □

**Theorem 3.5.** *Expression for elements of $_*{}^*$-inverse matrix has form*

$$
A^{-1_*{}^*} = \mathcal{H}\det(A,{}_*{}^*) \tag{3.13}
$$

*Proof.* (3.13) follows from (3.12). □

**Theorem 3.6.** *Expression for $\binom{a}{b}$-$_*{}^*$-quasideterminant can be evaluated by either form*[7]

$$
{}^j\det(A,{}_*{}^*)_i = {}^jA_i - {}^jA_{[i]}{}_*{}^*\left({}^{[j]}A_{[i]}\right)^{-1_*{}^*}{}_*{}^*{}^{[j]}A_i \tag{3.14}
$$

$$
{}^j\det(A,{}_*{}^*)_i = {}^jA_i - {}^jA_{[i]}{}_*{}^*\mathcal{H}\det\left({}^{[j]}A_{[i]},{}_*{}^*\right){}_*{}^*{}^{[j]}A_i \tag{3.15}
$$

*Proof.* Statement follows from (3.7) and (3.12). □

**Theorem 3.7.**

$$
{}_j\det\left(A^T,{}_*{}^*\right)^i = {}^j\det(A,{}^*{}_*)_i \tag{3.16}
$$

*Proof.* According to (3.12) and (1.2)

$$
{}_j\det\left(A^T,{}_*{}^*\right)^i = (\cdot{}^j_-((A^T)^{-1_*{}^*})\cdot{}^-_i)^{-1}
$$

Using theorem 2.12 we get

$$
{}_j\det\left(A^T,{}_*{}^*\right)^i = (\cdot{}^j_-((A^{-1^*{}_*})^T)\cdot{}^-_i)^{-1}
$$

---

[7] We can provide similar proof for $\binom{a}{b}$**-$^*{}_*$-quasideterminant**. However we can write corresponding statement using the duality principle. Thus, if we read equation (3.14) from right to left, we get equation

$$
{}^j\det(A,{}^*{}_*)_i = {}^jA_i - {}^{[j]}A_i{}^*{}_*\left({}^{[j]}A_{[i]}\right)^{-1^*{}_*}{}^*{}_*{}^jA_{[i]}
$$

$$
{}^j\det(A,{}^*{}_*)_i = {}^jA_i - {}^{[j]}A_i{}^*{}_*\mathcal{H}\det\left({}^{[j]}A_{[i]},{}^*{}_*\right){}^*{}_*{}^jA_{[i]}
$$

Using (1.1) we get

$$ {}_j\det\left(A^T,{}_*{}^*\right)^i = ({}^{\cdot}{}_j^{-}(A^{-1^*{}_*})^{i}_{\cdot -})^{-1} \tag{3.17}$$

Using (3.17), (1.2), (3.12) we get (3.16). □

The theorem 3.7 extends the duality principle stated in the theorem 2.14 to statements on quasideterminants and tells us that the same expression is ${}_*{}^*$-quasideterminant of matrix $A$ and ${}^*{}_*$-quasideterminant of matrix $A^T$. Using this theorem, we can write any statement for ${}^*{}_*$-matrix on the basis of similar statement for ${}_*{}^*$-matrix.

**Theorem 3.8** (duality principle)**.** *Let $\mathfrak{A}$ be true statement about matrix biring. If we exchange the same time*

- *${}_*$-row and ${}^*$-row*
- *${}_*{}^*$-quasideterminant and ${}^*{}_*$-quasideterminant*

*then we soon get true statement.*

**Theorem 3.9.**

$$ (mA)^{-1_*{}^*} = A^{-1_*{}^*}m^{-1} \tag{3.18}$$

$$ (Am)^{-1_*{}^*} = m^{-1}A^{-1_*{}^*} \tag{3.19}$$

*Proof.* To prove equation (3.18) we proceed by induction on size of the matrix.

Since

$$(mA)^{-1_*{}^*} = \left((mA)^{-1}\right) = \left(A^{-1}m^{-1}\right) = \left(A^{-1}\right)m^{-1} = A^{-1_*{}^*}m^{-1}$$

the statement is evident for $1\times 1$ matrix.

Let the statement holds for $(n-1)\times(n-1)$ matrix. Then from equation (3.1) it follows that

$$\begin{aligned}({}^I((mA)^{-1_*{}^*})_J)^{-1_*{}^*} &= {}^J(mA)_I - {}^J(mA)_{[I]}{}_*{}^*\left({}^{[J]}(mA)_{[I]}\right)^{-1_*{}^*}{}_*{}^*{}^{[J]}(mA)_I\\ &= m{}^JA_I - m\,{}^JA_{[I]}{}_*{}^*\left({}^{[J]}A_{[I]}\right)^{-1_*{}^*}\; m^{-1}{}_*{}^*m\,{}^{[J]}A_I\\ &= m\,{}^JA_I - m\,{}^JA_{[I]}{}_*{}^*\left({}^{[J]}A_{[I]}\right)^{-1_*{}^*}{}_*{}^*{}^{[J]}A_I\end{aligned}$$

$$ ({}^I((mA)^{-1_*{}^*})_J)^{-1_*{}^*} = m\,{}^I(A^{-1_*{}^*})_J \tag{3.20}$$

The equation (3.18) follows from the equation (3.20). In the same manner we prove the equation (3.19). □

**Theorem 3.10.** *Let*

$$ A = \begin{pmatrix}1 & 0\\ 0 & 1\end{pmatrix} \tag{3.21}$$

*Then*

$$ A^{-1_*{}^*} = \begin{pmatrix}1 & 0\\ 0 & 1\end{pmatrix} \tag{3.22}$$

$$ A^{-1^*{}_*} = \begin{pmatrix}1 & 0\\ 0 & 1\end{pmatrix} \tag{3.23}$$

*Proof.* It is clear from (3.8) and (3.11) that ${}^1(A^{-1_*{}^*})_1 = 1$ and ${}^2(A^{-1_*{}^*})_2 = 1$. However expression for ${}^2(A^{-1_*{}^*})_1$ and ${}^1(A^{-1_*{}^*})_2$ cannot be defined from (3.9) and (3.10) since ${}^2A_1 = {}^1A_2 = 0$. We can transform these expressions. For instance

$$\begin{aligned}
{}^2(A^{-1_*{}^*})_1 &= ({}^1A_2 - {}^1A_1({}^2A_1)^{-1}\ {}^2A_2)^{-1}\\
&= ({}^1A_1(({}^1A_1)^{-1}\ {}^1A_2 - ({}^2A_1)^{-1}\ {}^2A_2))^{-1}\\
&= (({}^2A_1)^{-1}\ {}^1A_1({}^2A_1({}^1A_1)^{-1}\ {}^1A_2 - {}^2A_2))^{-1}\\
&= ({}^1A_1({}^2A_1({}^1A_1)^{-1}\ {}^1A_2 - {}^2A_2))^{-1}\ {}^2A_1
\end{aligned}$$

It follows immediately that ${}^2(A^{-1_*{}^*})_1 = 0$. In the same manner we can find that ${}^1(A^{-1_*{}^*})_2 = 0$. This completes the proof of (3.22).

Equation (3.23) follows from (3.22), theorem 3.7 and symmetry of matrix (3.21). □

## 4. Reducible Biring

Let us study biring of matrices over field $F$. From the commutativity of product in the field it follows

$$A_*{}^*B = (A^c_a B^b_c) = (B^b_c A^c_a) = B^*{}_*A \tag{4.1}$$

**Definition 4.1. Reducible biring** is the biring which holds **condition of reducibility of products** (4.1). □

**Theorem 4.2.**

$$(A_*{}^*B)^T = B^T{}_*{}^*A^T \tag{4.2}$$

$$(A^*{}_*B)^T = B^{T*}{}_*A^T \tag{4.3}$$

*in the reducible biring.*

*Proof.* From (2.6) and (4.1) it follows that

$$(A_*{}^*B)^T = A^{T*}{}_*B^T = B^T{}_*{}^*A^T$$

We prove (4.3) the similar way. □

**Theorem 4.3.**

$$A^{-1^*{}_*} = A^{-1_*{}^*} \tag{4.4}$$

*in the reducible biring.*

*Proof.* From (2.22) and (4.1) it follows that

$$\delta = A_*{}^*A^{-1_*{}^*} = A^{-1_*{}^**}{}_*A \tag{4.5}$$

(4.4) follows from comparison of (4.5) and (2.23). □

## 5. Quasideterminant over Field

**Definition 5.1.** For matrices defined over a field we define a function which is called the **determinant of the matrix**,

$$\det() = 1 \tag{5.1}$$

$$\det a = \sum_a (-1)^{|a|+|b|} a_a^b \det a^{[a]}_{[b]} \tag{5.2}$$

□

**Theorem 5.2.**

$$\det (a,{}_*{}^*)^a_b = (-1)^{|a|+|b|} \frac{\det a}{\det a^{[b]}_{[a]}} \tag{5.3}$$

*for matrices defined over a field.*

*Proof.* We proceed by induction on $n$. Since the product in field is commutative we can reduce the expression (3.15)

$$\det (a,{}_*{}^*)^a_b = a^a_b - \frac{a^c_b a^a_e}{\det \left(a^{[b]}_{[a]},{}_*{}^*\right)^c_e} \tag{5.4}$$

For $n = 1$, we verify the statement of the theorem directly.

Let the statement of the theorem hold for $n = k - 1$. Then this statement holds for all minors of the $k \times k$ matrix $a$. We substitute equation (5.3) in equation (5.4)

$$\begin{aligned}
\det (a,{}_*{}^*)^a_b = a^a_b &- \sum_{|c|<|a|} \sum_{|e|<|b|} \frac{a^c_b a^a_e}{(-1)^{|c|+|e|} \frac{det\ a^{[b]}_{[a]}}{det\ a^{[b,e]}_{[a,c]}}} \\
&+ \sum_{|c|>|a|} \sum_{|e|<|b|} \frac{a^c_b a^a_e}{(-1)^{|c|+|e|} \frac{det\ a^{[b]}_{[a]}}{det\ a^{[b,e]}_{[a,c]}}} \\
&+ \sum_{|c|<|a|} \sum_{|e|>|b|} \frac{a^c_b a^a_e}{(-1)^{|c|+|e|} \frac{det\ a^{[b]}_{[a]}}{det\ a^{[b,e]}_{[a,c]}}} \\
&- \sum_{|c|>|a|} \sum_{|e|>|b|} \frac{a^c_b a^a_e}{(-1)^{|c|+|e|} \frac{det\ a^{[b]}_{[a]}}{det\ a^{[b,e]}_{[a,c]}}}
\end{aligned}$$

$$\begin{aligned}
\det (a,{}_*{}^*)^a_b = a^a_b \quad &- \sum_{|c|<|a|} \sum_{|e|<|b|} \frac{(-1)^{|a|+|b|+|e|+|c|} a^c_b a^a_e det\ a^{[b,e]}_{[a,c]}}{(-1)^{|a|+|b|} det\ a^{[b]}_{[a]}} \\
&+ \sum_{|c|>|a|} \sum_{|e|<|b|} \frac{(-1)^{|a|+|b|+|e|+|c|} a^c_b a^a_e det\ a^{[b,e]}_{[a,c]}}{(-1)^{|a|+|b|} det\ a^{[b]}_{[a]}} \\
&+ \sum_{|c|<|a|} \sum_{|e|>|b|} \frac{(-1)^{|a|+|b|+|e|+|c|} a^c_b a^a_e det\ a^{[b,e]}_{[a,c]}}{(-1)^{|a|+|b|} det\ a^{[b]}_{[a]}} \\
&- \sum_{|c|>|a|} \sum_{|e|>|b|} \frac{(-1)^{|a|+|b|+|e|+|c|} a^c_b a^a_e det\ a^{[b,e]}_{[a,c]}}{(-1)^{|a|+|b|} det\ a^{[b]}_{[a]}}
\end{aligned} \tag{5.5}$$

Substituting equation

$$\det a^{[b]}_{[c]} = \begin{cases} \left(\sum_{|e|<|b|} - \sum_{|e|>|b|}\right)(-1)^{|a|+|e|} a^a_e det\ a^{[b,e]}_{[a,c]} & c > a \\ \left(\sum_{|e|>|b|} - \sum_{|e|<|b|}\right)(-1)^{|a|+|e|} a^a_e det\ a^{[b,e]}_{[a,c]} & c < a \end{cases}$$

into the equation (5.5), we get

$$\det (a,{}_*{}^*)^a_b = a^a_b + \left(\sum_{|c|<|a|} + \sum_{|c|>|a|}\right) \frac{(-1)^{|b|+|c|} a^c_b det\ a^{[b]}_{[c]}}{(-1)^{|a|+|b|} det\ a^{[b]}_{[a]}}$$

Adding fractions, we get

$$\det (a,{}_*{}^*)^a_b = (-1)^{|a|+|b|} \sum_{c\in M} \frac{(-1)^{|b|+|c|} a^c_b \det a^{[b]}_{[c]}}{det\ a^{[b]}_{[a]}} \tag{5.6}$$

Substituting equation (5.2) into (5.6), we prove that the statement of the theorem hold for $n = k$. □

# 7. Index

# 8. Special Symbols and Notations

# Бикольцо матриц

Александр Клейн

Аннотация. Матрицы допускают две операции произведения, связанные операцией транспонирования. Бикольцо - это алгебра, определяющая на множестве две взаимносвязанные структуры кольца. Согласно каждому виду произведения мы можем расширить определение квазидетерминанта, данное в [1, 2], и определить два разных вида квазидетерминанта.

## Содержание

## 1. Концепция обобщённого индекса

Изучая тензорное исчисление, мы начинаем с изучения одновалентных ковариантного и контравариантного тензоров. Несмотря на различие свойств, оба эти объекта являются элементами соответствующих векторных пространств. Если мы введём обобщённый индекс по правилу $a^{i} = a^{i}$, $b^{i} = b^{\cdot}{}^{-}_{i}$, то мы видим, что эти тензоры ведут себя одинаково. Например, преобразование ковариантного тензора принимает форму

$$b'^{i} = b'^{\cdot}{}^{-}_{i} = f^{\cdot}{}^{-}_{i}{}^{\cdot}{}^{j}_{-} b^{\cdot}{}^{-}_{j} = f^{i}_{j} b^{j}$$

Это сходство идёт сколь угодно далеко, так как тензоры также порождают векторное пространство.

Эти наблюдения сходства свойств ковариантного и контравариантного тензоров приводят нас к концепции обобщённого индекса.[1] Я пользуюсь символом · перед обобщённым индексом, когда мне необходимо описать его структуру.

Aleks_Kleyn@MailAPS.org.
http://sites.google.com/site/AleksKleyn/.
http://arxiv.org/a/kleyn_a_1.
http://AleksKleyn.blogspot.com/.

[1]Набор команд для работы с обобщёнными индексами, также как и другие команды, используемые в этой и последующих статьях, можно найти в файле Commands.tex.

Я помещаю символ $'-'$ на месте индекса, позиция которого изменилась. Например, если исходное выражение было $a_{ij}$, я пользуюсь записью $a_{i-}^{\ j}$ вместо записи $a_i^j$.

Хотя структура обобщённого индекса произвольна, мы будем предполагать, что существует взаимно однозначное отображение отрезка натуральных чисел 1, ..., $n$ на множество значений индекса. Пусть $i$ - множество значений индекса $i$. Мы будем обозначать мощность этого множества символом $|i|$ и будем полагать $|i| = n$. Если нам надо перечислить элементы $a_i$, мы будем пользоваться обозначением $a_1$, ..., $a_n$.

Представление координат вектора в форме матрицы позволяет сделать запись более компактной. Вопрос о представлении вектора как строка или столбец матрицы является вопросом соглашения. Мы можем распространить концепцию обобщённого индекса на элементы матрицы. Матрица - это двумерная таблица, строки и столбцы которой занумерованы обобщёнными индексами. Для представления матрицы мы будем пользоваться одним из следующих представлений:

**Стандартное представление:** в этом случае мы представляем элементы матрицы $A$ в виде $A_b^a$.

**Альтернативное представление:** в этом случае мы представляем элементы матрицы $A$ в виде ${}^aA_b$ или ${}_bA^a$.

Так как мы пользуемся обобщёнными индексами, мы не можем сказать, нумерует ли индекс $a$ строки матрицы или столбцы, до тех пор, пока мы не знаем структуры индекса.

Мы могли бы пользоваться терминами $*$-столбец и $*$-строка, которые более близки традиционным. Однако как мы увидим ниже для нас несущественна форма представления матрицы. Для того, чтобы обозначения, предлагаемые ниже, были согласованы с традиционными, мы будем предполагать, что матрица представлена в виде

$$A = \begin{pmatrix} {}^1A_1 & ... & {}^1A_n \\ ... & ... & ... \\ {}^mA_1 & ... & mA_n \end{pmatrix}$$

**Определение 1.1.** *Я использую следующие имена и обозначения различных миноров матрицы $A$*

$A_a$ **:** $^*$**-строка** *с индексом $a$ является обобщением столбца матрицы. Верхний индекс перечисляет элементы $^*$-строки, нижний индекс перечисляет $^*$-строки.*

$A_T$ **:** *минор, полученный из $A$ выбором $^*$-строк с индексом из множества $T$*

$A_{[a]}$ **:** *минор, полученный из $A$ удалением $^*$-строки $A_a$*

$A_{[T]}$ **:** *минор, полученный из $A$ удалением $^*$-строк с индексом из множества $T$*

${}^bA$ **:** $_*$**-строка** *с индексом $b$ является обобщением строки матрицы. Нижний индекс перечисляет элементы $_*$-строки, верхний индекс перечисляет $_*$-строки.*

${}^SA$ **:** *минор, полученный из $A$ выбором $_*$-строк с индексом из множества $S$*

${}^{[b]}A$ **:** *минор, полученный из $A$ удалением $_*$-строки ${}^bA$*

$^{[S]}A$ **:** минор, полученный из $A$ удалением $_*$-строк с индексом из множества $S$

□

*Замечание* 1.2. Мы будем комбинировать запись индексов. Так ${}^bA_a$ является $1\times 1$ минором. Одновременно, это обозначение элемента матрицы. Это позволяет отождествить $1\times 1$ матрицу и её элемент. Индекс $a$ является номером $^*$-строки матрицы и индекс $b$ является номером $_*$-строки матрицы. □

Каждая форма записи матрицы имеет свои преимущества. Стандартная форма более естественна, когда мы изучаем теорию матриц. Альтернативная форма записи делает выражения более ясными в теории векторных пространств. Распространив альтернативную запись индексов на произвольные тензоры, мы сможем лучше понять взаимодействие различных геометрических объектов. Опираясь на принцип двойственности (теорема 2.14), мы можем расширить наши выразительные возможности.

*Замечание* 1.3. Мы можем договориться, что при чтении мы произносим символ $^*$- как $c$- и символ $_*$- как $r$-, формируя тем самым названия **$c$-строка** и **$r$-строка**. В последующем мы распространим это соглашение на другие элементы линейной алгебры. Я буду пользоваться этим соглашением при составлении индекса. □

Так как транспонирование матрицы меняет местами $_*$-строки и $^*$-строки, то мы получаем равенство

$$ {}_i(A^T)^j = {}^iA_j \tag{1.1} $$

*Замечание* 1.4. Как видно из равенства (1.1), для нас несущественно, с какой стороны мы указываем номер $_*$-строки и с какой стороны мы указываем номер $^*$-строки. Это связано с тем, что мы можем нумеровать элементы матрицы различными способами. Если мы хотим указывать номера $_*$-строки и $^*$-строки согласно определению 1.1, то равенство (1.1) примет вид

$$ {}^j(A^T)_i = {}^iA_j $$

В стандартном представлении равенство (1.1) примет вид

$$ (A^T)_i^j = A_j^i $$

□

Мы называем матрицу[2]

$$ \mathcal{H}A = ({}^j\mathcal{H}A_i) = ((\cdot{}_i^{-}A\cdot{}_{-}^{j})^{-1}) \tag{1.2} $$

**обращением Адамара матрицы** $A = ({}_bA^a)$ ([2]-page 4).

Я пользуюсь Эйнштейновским соглашением о суммировании. Это означает, что, когда индекс присутствует в выражении дважды и множество индексов

[2] Запись $(\cdot{}_i^{-}A\cdot{}_{-}^{j})^{-1}$ означает, что при обращении Адамара столбцы и строки меняются местами. Мы можем формально записать это выражение следующим образом

$$ (\cdot{}_i^{-}A\cdot{}_{-}^{j})^{-1} = \frac{1}{{}^iA_j} $$

известно, у меня есть сумма по этому индексу. Я буду явно указывать множество индексов, если это необходимо. Кроме того, в этой статье я использую туже корневую букву для матрицы и её элементов.

Мы будем изучать матрицы, элементы которых принадлежат телу $D$. Мы также будем иметь в виду, что вместо тела $D$ мы можем в тексте писать поле $F$. Мы будем явно писать поле $F$ в тех случаях, когда коммутативность будет порождать новые детали. Обозначим через 1 единичный элемент тела $D$.

Пусть $I$, $|I| = n$ - множество индексов. **Символ Кронекера** определён равенством

$$\delta^i_j = \begin{cases} 1 & i = j \\ 0 & i \neq j \end{cases} \qquad i, j \in I \tag{1.3}$$

## 2. Бикольцо

Мы будем рассматривать матрицы, элементы которых принадлежат телу $D$.

Произведение матриц связано с произведением гомоморфизмов векторных пространств над полем. Согласно традиции произведение матриц $A$ и $B$ определено как произведение $_*$-строк матрицы $A$ и $^*$-строк матрицы $B$. Условность этого определения становится очевидной, если мы обратим внимание, что $_*$-строка матрицы $A$ может быть столбцом этой матрицы. В этом случае мы умножаем столбцы матрицы $A$ на строки матрицы $B$. Таким образом, мы можем определить два вида произведения матриц. Чтобы различать эти произведения, мы вводим новые обозначения.[3]

**Определение 2.1.** **$_*{}^*$-произведение матриц** $A$ и $B$ имеет вид

$$\begin{cases} A_*{}^*B = ({}^aA_c\ {}^cB_b) \\ {}^a(A_*{}^*B)_b = {}^aA_c\ {}^cB_b \end{cases} \tag{2.1}$$

и может быть выражено как произведение $_*$-строки матрицы $A$ и $^*$-строки матрицы $B$.[4] □

**Определение 2.2.** **$^*{}_*$-произведение матриц** $A$ и $B$ имеет вид

$$\begin{cases} A^*{}_*B = ({}_aA^c\ {}_cB^b) \\ {}_a(A^*{}_*B)^b = {}_aA^c\ {}_cB^b \end{cases} \tag{2.2}$$

и может быть выражено как произведение $^*$-строки матрицы $A$ на $_*$-строку матрицы $B$.[5] □

---

[3] Для совместимости обозначений с существующими мы будем иметь в виду $_*{}^*$-произведение, когда нет явных обозначений.

[4] В альтернативной форме операция состоит из двух символов $*$, которые записываются на месте индекса суммирования. В стандартной форме операция имеет вид

$$\begin{cases} A_*{}^*B = (A^a_c B^c_b) \\ (A_*{}^*B)^a_b = A^a_c B^c_b \end{cases}$$

и может быть интерпретирована как символическая запись

$$A_*{}^*B = A_*B^*$$

где мы записываем символ $*$ на месте индекса, по которому предполагается суммирование.

[5] В альтернативной форме операция состоит из двух символов $*$, которые записываются на месте индекса суммирования. В стандартной форме операция имеет вид

$$\begin{cases} A^*{}_*B = (A^c_a\ B^b_c) \\ (A^*{}_*B)^a_b = A^c_a\ B^b_c \end{cases}$$

*Замечание* 2.3. Мы будем пользоваться символом ${}_*{}^*$- или ${}^*{}_*$- в имени свойств каждого произведения и в обозначениях. Согласно замечанию 1.3 мы можем читать символы ${}_*{}^*$ и ${}^*{}_*$ как $rc$-произведение и $cr$-произведение. Это правило мы распространим на последующую терминологию. □

*Замечание* 2.4. Также как и в замечании 1.4, я хочу обратить внимание на то, что я меняю нумерацию элементов матрицы. Если мы хотим указывать номера ${}_*$-строки и ${}^*$-строки согласно определению 1.1, то равенство (2.2) примет вид

$$
{}^b(A^*{}_*B)_a = {}^cA_a\ {}^bB_c \tag{2.3}
$$

Однако формат равенства (2.3) несколько необычен. □

Множество $n \times n$ матриц замкнуто относительно ${}_*{}^*$-произведения и ${}^*{}_*$-произведения, а также относительно суммы, определённой согласно правилу

$$
(A+B)_a^b = A_a^b + B_a^b
$$

**Теорема 2.5.**

$$
(A_*{}^*B)^T = A^T{}^*{}_*B^T \tag{2.4}
$$

*Доказательство.* Цепочка равенств

$$
\begin{aligned}
{}_a((A_*{}^*B)^T)^b &= {}^a(A_*{}^*B)_b \\
&= {}^aA_c\ {}^cB_b \\
&= {}_a(A^T)^c{}_c(B^T)^b \\
&= {}_a((A^T){}^*{}_*(B^T))^b
\end{aligned} \tag{2.5}
$$

следует из (1.1), (2.1) и (2.2). Равенство (2.4) следует из (2.5). □

Матрица $\delta = (\delta_a^c)$ является единицей для обоих произведений.

**Определение 2.6. Бикольцо** $\mathcal{A}$ - это множество, на котором мы определили унарную операцию, называемую транспозицией, и три бинарных операции, называемые ${}_*{}^*$-произведение, ${}^*{}_*$-произведение и сумма, такие что

- ${}_*{}^*$-произведение и сумма определяют структуру кольца на $\mathcal{A}$
- ${}^*{}_*$-произведение и сумма определяют структуру кольца на $\mathcal{A}$
- оба произведения имеют общую единицу $\delta$
- произведения удовлетворяют равенству

$$
(A_*{}^*B)^T = A^T{}^*{}_*B^T \tag{2.6}
$$

- транспозиция единицы есть единица

$$
\delta^T = \delta \tag{2.7}
$$

- двойная транспозиция есть исходный элемент

$$
(A^T)^T = A \tag{2.8}
$$

□

---

и может быть интерпретирована как символическую запись

$$
A^*{}_*B = A^*B_*
$$

где мы записываем символ $*$ на месте индекса, по которому предполагается суммирование.

**Теорема 2.7.**

$$(A^*{}_*B)^T = (A^T)_*{}^*(B^T) \tag{2.9}$$

*Доказательство.* Мы можем доказать (2.9) в случае матриц тем же образом, что мы доказали (2.6). Тем не менее для нас более важно показать, что (2.9) следует непосредственно из (2.6).

Применяя (2.8) к каждому слагаемому в левой части (2.9), мы получим

$$(A^*{}_*B)^T = ((A^T)^T{}^*{}_*(B^T)^T)^T \tag{2.10}$$

Из (2.10) и (2.6) следует, что

$$(A^*{}_*B)^T = ((A^T{}_*{}^*B^T)^T)^T \tag{2.11}$$

(2.9) следует из (2.11) и (2.8). □

**Определение 2.8.** Мы определим ${}_*{}^*$**-степень** элемента $A$ бикольца $\mathcal{A}$, пользуясь рекурсивным правилом

$$A^{0_*{}^*} = \delta \tag{2.12}$$

$$A^{n_*{}^*} = A^{n-1_*{}^*}{}_*{}^*A \tag{2.13}$$

□

**Определение 2.9.** Мы определим ${}^*{}_*$**-степень** элемента $A$ бикольца $\mathcal{A}$, пользуясь рекурсивным правилом

$$A^{0^*{}_*} = \delta \tag{2.14}$$

$$A^{n^*{}_*} = A^{n-1^*{}_*}{}^*{}_*A \tag{2.15}$$

□

**Теорема 2.10.**

$$(A^T)^{n_*{}^*} = (A^{n^*{}_*})^T \tag{2.16}$$

$$(A^T)^{n^*{}_*} = (A^{n_*{}^*})^T \tag{2.17}$$

*Доказательство.* Мы проведём доказательство индукцией по $n$.

При $n = 0$ утверждение непосредственно следует из равенств (2.12), (2.14) и (2.7).

Допустим утверждение справедливо при $n = k - 1$

$$(A^T)^{n-1_*{}^*} = (A^{n-1^*{}_*})^T \tag{2.18}$$

Из (2.13) следует

$$(A^T)^{k_*{}^*} = (A^T)^{k-1_*{}^*}{}_*{}^*A^T \tag{2.19}$$

Из (2.19) и (2.18) следует

$$(A^T)^{k_*{}^*} = (A^{k-1^*{}_*})^T{}_*{}^*A^T \tag{2.20}$$

Из (2.20) и (2.9) следует

$$(A^T)^{k_*{}^*} = (A^{k-1^*{}_*}{}^*{}_*A)^T \tag{2.21}$$

Из (2.19) и (2.15) следует (2.16).

Мы можем доказать (2.17) подобным образом. □

**Определение 2.11.** Элемент $A^{-1_*{}^*}$ бикольца $\mathcal{A}$ - это **${}_*{}^*$-обратный элемент** элемента $A$, если

$$A_*{}^*A^{-1_*{}^*}=\delta \tag{2.22}$$

Элемент $A^{-1^*{}_*}$ бикольца $\mathcal{A}$ - это **${}^*{}_*$-обратный элемент** элемента $A$, если

$$A^*{}_*A^{-1^*{}_*}=\delta \tag{2.23}$$

□

**Теорема 2.12.** *Предположим, что элемент $A\in\mathcal{A}$ имеет ${}_*{}^*$-обратный элемент. Тогда транспонированный элемент $A^T$ имеет ${}^*{}_*$-обратный элемент и эти элементы удовлетворяют равенству*

$$(A^T)^{-1^*{}_*}=(A^{-1_*{}^*})^T \tag{2.24}$$

*Предположим, что элемент $A\in\mathcal{A}$ имеет ${}^*{}_*$-обратный элемент. Тогда транспонированный элемент $A^T$ имеет ${}_*{}^*$-обратный элемент и эти элементы удовлетворяет равенству*

$$(A^T)^{-1_*{}^*}=(A^{-1^*{}_*})^T \tag{2.25}$$

*Доказательство.* Если мы возьмём транспонирование обеих частей (2.22) и применим (2.7), мы получим

$$(A_*{}^*A^{-1_*{}^*})^T=\delta^T=\delta$$

Применяя (2.6), мы получим

$$\delta=A^T{}^*{}_*(A^{-1_*{}^*})^T \tag{2.26}$$

(2.24) следует из сравнения (2.23) и (2.26).

Мы можем доказать (2.25) подобным образом. □

Теоремы 2.5, 2.7, 2.10 и 2.12 показывают, что существует двойственность между ${}_*{}^*$-произведением и ${}^*{}_*$-произведением. Мы можем объединить эти утверждения.

**Теорема 2.13 (принцип двойственности для бикольца).** *Пусть $\mathfrak{A}$ - истинное утверждение о бикольце $\mathcal{A}$. Если мы заменим одновременно*

- *$A\in\mathcal{A}$ и $A^T$*
- *${}_*{}^*$-произведение и ${}^*{}_*$-произведение*

*то мы снова получим истинное утверждение.*

**Теорема 2.14 (принцип двойственности для бикольца матриц).** *Пусть $\mathcal{A}$ является бикольцом матриц. Пусть $\mathfrak{A}$ - истинное утверждение о матрицах. Если мы заменим одновременно*

- *${}^*$-строки и ${}_*$-строки всех матриц*
- *${}_*{}^*$-произведение и ${}^*{}_*$-произведение*

*то мы снова получим истинное утверждение.*

*Доказательство.* Непосредственное следствие теоремы 2.13. □

*Замечание* 2.15. В выражении

$$A_*{}^*B_*{}^*C$$

мы выполняем операцию умножения слева направо. Однако мы можем выполнять операцию умножения справа налево. В традиционной записи это выражение примет вид

$$C^*{}_*B^*{}_*A$$

Мы сохраним правило, что показатель степени записывается справа от выражения. Если мы пользуемся стандартным представлением, то индексы также записываются справа от выражения. Если мы пользуемся альтернативным представлением, то индексы читаются в том же порядке, что и символы операции и корневые буквы. Например, если исходное выражение имеет вид

$$A^{-1_*{}^*}{}_*{}^*B_a$$

то выражение, читаемое справа налево, примет вид

$$B_a{}^*{}_*A^{-1^*{}_*}$$

в стандартном представлении либо примет вид

$${}_aB^*{}_*A^{-1^*{}_*}$$

в альтернативном представлении.

Если задать порядок, в котором мы записываем индексы, то мы можем утверждать, что мы читаем выражение сверху вниз, читая сперва верхние индексы, потом нижние. Договорившись, что это стандартная форма чтения, мы можем прочесть выражение снизу вверх. При этом мы дополним правило, что символы операции также читаются в том же направлении, что и индексы. Например, выражение

$$A^a{}_*{}^*B^{-1_*{}^*} = C^a$$

прочтённое снизу вверх, в стандартной форме имеет вид

$$A_a{}^*{}_*B^{-1^*{}_*} = C_a$$

Согласно принципу двойственности, если верно одно утверждение, то верно и другое. □

**Теорема 2.16.** *Если матрица $A$ имеет ${}_*{}^*$-обратную матрицу, то для любых матриц $B$ и $C$ из равенства*

$$B_*{}^*A = C_*{}^*A \tag{2.27}$$

*следует равенство*

$$B = C \tag{2.28}$$

*Доказательство.* Равенство (2.28) следует из (2.27), если обе части равенства (2.27) умножить на $A^{-1_*{}^*}$. □

## 3. Квазидетерминант

**Теорема 3.1.** *Предположим, что $n\times n$ матрица $A$ имеет ${}_*{}^*$-обратную матрицу.*[6] *Тогда $k\times k$ минор ${}_*{}^*$-обратной матрицы удовлетворяет*

$$\left({}^I(A^{-1_*{}^*})_J\right)^{-1_*{}^*}={}^JA_I-{}^JA_{[I]}{}_*{}^*\left({}^{[J]}A_{[I]}\right)^{-1_*{}^*}{}_*{}^*{}^{[J]}A_I \tag{3.1}$$

*Доказательство.* Определение (2.22) ${}_*{}^*$-обратной матрицы приводит к системе линейных уравнений

$$ {}^{[J]}A_{[I]}{}_*{}^*{}^{[I]}(A^{-1_*{}^*})_J+{}^{[J]}A_I{}_*{}^*{}^I(A^{-1_*{}^*})_J=0 \tag{3.2}$$

$$ {}^JA_{[I]}{}_*{}^*{}^{[I]}(A^{-1_*{}^*})_J+{}^JA_I{}_*{}^*{}^I(A^{-1_*{}^*})_J=\delta \tag{3.3}$$

Мы умножим (3.2) на $\left({}^{[J]}A_{[I]}\right)^{-1_*{}^*}$

$$ {}^{[I]}(A^{-1_*{}^*})_J+\left({}^{[J]}A_{[I]}\right)^{-1_*{}^*}{}_*{}^*{}^{[J]}A_I{}_*{}^*{}^I(A^{-1_*{}^*})_J=0 \tag{3.4}$$

Теперь мы можем подставить (3.4) в (3.3)

$$ -{}^JA_{[I]}{}_*{}^*\left({}^{[J]}A_{[I]}\right)^{-1_*{}^*}{}_*{}^*{}^{[J]}A_I{}_*{}^*{}^I(A^{-1_*{}^*})_J+{}^JA_I{}_*{}^*{}^I(A^{-1_*{}^*})_J=\delta \tag{3.5}$$

(3.1) следует из (3.5). □

**Следствие 3.2.** *Предположим, что $n\times n$ матрица $A$ имеет ${}_*{}^*$-обратную матрицу. Тогда элементы ${}_*{}^*$-обратной матрицы удовлетворяют равенству*[2]

$$ {}^i(A^{-1_*{}^*})_j=\left({}^jA_i-{}^jA_{[i]}{}_*{}^*\left({}^{[j]}A_{[i]}\right)^{-1_*{}^*}{}_*{}^*{}^{[j]}A_i\right)^{-1} \tag{3.6}$$

$$ {}^j\left(\mathcal{H}A^{-1_*{}^*}\right)_i={}^jA_i-{}^jA_{[i]}{}_*{}^*\left({}^{[j]}A_{[i]}\right)^{-1_*{}^*}{}_*{}^*{}^{[j]}A_i \tag{3.7}$$

□

**Пример 3.3.** Рассмотрим матрицу

$$\begin{pmatrix}{}^1A_1 & {}^1A_2\\ {}^2A_1 & {}^2A_2\end{pmatrix}$$

Согласно (3.6)

$$ {}^1(A^{-1_*{}^*})_1=({}^1A_1-{}^1A_2({}^2A_2)^{-1}\,{}^2A_1)^{-1} \tag{3.8}$$

$$ {}^2(A^{-1_*{}^*})_1=({}^1A_2-{}^1A_1({}^2A_1)^{-1}\,{}^2A_2)^{-1} \tag{3.9}$$

$$ {}^1(A^{-1_*{}^*})_2=({}^2A_1-{}^2A_2({}^1A_2)^{-1}\,{}^1A_1)^{-1} \tag{3.10}$$

$$ {}^2(A^{-1_*{}^*})_2=({}^2A_2-{}^2A_1({}^1A_1)^{-1}\,{}^1A_2)^{-1} \tag{3.11}$$

Рассмотрим произведение матриц

$$C=\begin{pmatrix}({}^1A_1-{}^1A_2({}^2A_2)^{-1}\,{}^2A_1)^{-1} & ({}^2A_1-{}^2A_2({}^1A_2)^{-1}\,{}^1A_1)^{-1}\\ ({}^1A_2-{}^1A_1({}^2A_1)^{-1}\,{}^2A_2)^{-1} & ({}^2A_2-{}^2A_1({}^1A_1)^{-1}\,{}^1A_2)^{-1}\end{pmatrix}{}_*{}^*\begin{pmatrix}{}^1A_1 & {}^1A_2\\ {}^2A_1 & {}^2A_2\end{pmatrix}$$

[6]Это утверждение и его доказательство основаны на утверждении 1.2.1 из [1] (page 8) для матриц над свободным кольцом с делением.

Из непосредственных вычислений следует

$$
\begin{aligned}
{}^1C_1 &= ({}^1A_1 - {}^1A_2({}^2A_2)^{-1}\ {}^2A_1)^{-1}\ {}^1A_1 + ({}^2A_1 - {}^2A_2({}^1A_2)^{-1}\ {}^1A_1)^{-1}\ {}^2A_1\\
&= ({}^1A_2(({}^1A_2)^{-1}\ {}^1A_1 - ({}^2A_2)^{-1}\ {}^2A_1))^{-1}\ {}^1A_1\\
&+ ({}^2A_2(({}^2A_2)^{-1}\ {}^2A_1 - ({}^1A_2)^{-1}\ {}^1A_1))^{-1}\ {}^2A_1\\
&= (({}^1A_2)^{-1}\ {}^1A_1 - ({}^2A_2)^{-1}\ {}^2A_1)^{-1}\ ({}^1A_2)^{-1}\ {}^1A_1\\
&+ (({}^2A_2)^{-1}\ {}^2A_1 - ({}^1A_2)^{-1}\ {}^1A_1)^{-1}\ ({}^2A_2)^{-1}\ {}^2A_1\\
&= 1
\end{aligned}
$$

$$
\begin{aligned}
{}^1C_2 &= ({}^1A_1 - {}^1A_2({}^2A_2)^{-1}\ {}^2A_1)^{-1}\ {}^1A_2 + ({}^2A_1 - {}^2A_2({}^1A_2)^{-1}\ {}^1A_1)^{-1}\ {}^2A_2\\
&= ({}^1A_2(({}^1A_2)^{-1}\ {}^1A_1 - ({}^2A_2)^{-1}\ {}^2A_1))^{-1}\ {}^1A_2\\
&+ ({}^2A_2(({}^2A_2)^{-1}\ {}^2A_1 - ({}^1A_2)^{-1}\ {}^1A_1))^{-1}\ {}^2A_2\\
&= (({}^1A_2)^{-1}\ {}^1A_1 - ({}^2A_2)^{-1}\ {}^2A_1)^{-1}\ ({}^1A_2)^{-1}\ {}^1A_2\\
&+ (({}^2A_2)^{-1}\ {}^2A_1 - ({}^1A_2)^{-1}\ {}^1A_1)^{-1}\ ({}^2A_2)^{-1}\ {}^2A_2\\
&= 0
\end{aligned}
$$

$$
\begin{aligned}
{}^2C_1 &= ({}^1A_2 - {}^1A_1({}^2A_1)^{-1}\ {}^2A_2)^{-1}\ {}^1A_1 + ({}^2A_2 - {}^2A_1({}^1A_1)^{-1}\ {}^1A_2)^{-1}\ {}^2A_1\\
&= ({}^1A_1(({}^1A_1)^{-1}\ {}^1A_2 - ({}^2A_1)^{-1}\ {}^2A_2)^{-1}\ {}^1A_1\\
&+ ({}^2A_1(({}^2A_1)^{-1}\ {}^2A_2 - ({}^1A_1)^{-1}\ {}^1A_2))^{-1}\ {}^2A_1\\
&= (({}^1A_1)^{-1}\ {}^1A_2 - ({}^2A_1)^{-1}\ {}^2A_2)^{-1}\ ({}^1A_1)^{-1}\ {}^1A_1\\
&+ (({}^2A_1)^{-1}\ {}^2A_2 - ({}^1A_1)^{-1}\ {}^1A_2)^{-1}\ ({}^2A_1)^{-1}\ {}^2A_1\\
&= 0
\end{aligned}
$$

$$
\begin{aligned}
{}^2C_2 &= ({}^1A_2 - {}^1A_1({}^2A_1)^{-1}\ {}^2A_2)^{-1}\ {}^1A_2 + ({}^2A_2 - {}^2A_1({}^1A_1)^{-1}\ {}^1A_2)^{-1}\ {}^2A_2\\
&= ({}^1A_1(({}^1A_1)^{-1}\ {}^1A_2 - ({}^2A_1)^{-1}\ {}^2A_2)^{-1}\ {}^1A_2\\
&+ ({}^2A_1(({}^2A_1)^{-1}\ {}^2A_2 - ({}^1A_1)^{-1}\ {}^1A_2))^{-1}\ {}^2A_2\\
&= (({}^1A_1)^{-1}\ {}^1A_2 - ({}^2A_1)^{-1}\ {}^2A_2)^{-1}\ ({}^1A_1)^{-1}\ {}^1A_2\\
&+ (({}^2A_1)^{-1}\ {}^2A_2 - ({}^1A_1)^{-1}\ {}^1A_2)^{-1}\ ({}^2A_1)^{-1}\ {}^2A_2\\
&= 1
\end{aligned}
$$

□

Согласно [1], page 3 у нас нет определения детерминанта в случае тела. Тем не менее, мы можем определить квазидетерминант, который в конечном итоге даёт похожую картину. В определении, данном ниже, мы следуем определению [1]-1.2.2.

**Определение 3.4.** $\binom{j}{i}$**-**${}_*{}^*$**-квазидетерминант** $n \times n$ матрицы $A$ - это формальное выражение[2]

$$
{}^j\det(A, {}_*{}^*)_i = {}^j\left(\mathcal{H}A^{-1_*{}^*}\right)_i \tag{3.12}
$$

Согласно замечанию 1.2 мы можем рассматривать $\binom{j}{i}$-${}_*{}^*$-квазидетерминант как элемент матрицы $\det(a, {}_*{}^*)$, которую мы будем называть ${}_*{}^*$**-квазидетерминантом**. □

**Теорема 3.5.** *Выражение для элементов ${}_*{}^*$-обратной матрицы имеет вид*

$$A^{-1_*{}^*}=\mathcal{H}\det\left(A,{}_*{}^*\right)\tag{3.13}$$

*Доказательство.* (3.13) следует из (3.12). □

**Теорема 3.6.** *Выражение для $\binom{a}{b}$-${}_*{}^*$-квазидетерминант имеет любую из следующих форм*[7]

$$^j\det\left(A,{}_*{}^*\right)_i={}^jA_i-{}^jA_{[i]}{}_*{}^*\left({}^{[j]}A_{[i]}\right)^{-1_*{}^*}{}_*{}^*{}^{[j]}A_i\tag{3.14}$$

$$^j\det\left(A,{}_*{}^*\right)_i={}^jA_i-{}^jA_{[i]}{}_*{}^*\mathcal{H}\det\left({}^{[j]}A_{[i]},{}_*{}^*\right){}_*{}^*{}^{[j]}A_i\tag{3.15}$$

*Доказательство.* Утверждение следует из (3.7) и (3.12). □

**Теорема 3.7.**

$$_j\det\left(A^T,{}_*{}^*\right)^i={}^j\det\left(A,{}^*{}_*\right)_i\tag{3.16}$$

*Доказательство.* Согласно (3.12) и (1.2)

$$_j\det\left(A^T,{}_*{}^*\right)^i=\left({}^{\,j}_{\cdot\,-}\left(\left(A^T\right)^{-1_*{}^*}\right)^{\cdot\,-}_{\;\;i}\right)^{-1}$$

Пользуясь теоремой 2.12, мы получим

$$_j\det\left(A^T,{}_*{}^*\right)^i=\left({}^{\,j}_{\cdot\,-}\left(\left(A^{-1^*{}_*}\right)^T\right)^{\cdot\,-}_{\;\;i}\right)^{-1}$$

Пользуясь (1.1), мы имеем

$$_j\det\left(A^T,{}_*{}^*\right)^i=\left({}^{\cdot\,-}_{\;\;j}\left(A^{-1^*{}_*}\right){}^{\,i}_{\cdot\,-}\right)^{-1}\tag{3.17}$$

Пользуясь (3.17), (1.2), (3.12), мы получим (3.16). □

Теорема 3.7 расширяет принцип двойственности, теорема 2.14, на утверждения о квазидетерминантах и утверждает, что одно и тоже выражение является ${}_*{}^*$-квазидетерминантом матрицы $A$ и ${}^*{}_*$-квазидетерминантом матрицы $A^T$. Пользуясь этой теоремой, мы можем записать любое утверждение о ${}^*{}_*$-матрице, опираясь на подобное утверждение о ${}_*{}^*$-матрице.

**Теорема 3.8** (принцип двойственности)**.** *Пусть $\mathfrak{A}$ - истинное утверждение о бикольце матриц. Если мы одновременно заменим*

- *${}_*$-строку и ${}^*$-строку*
- *${}_*{}^*$-квазидетерминант и ${}^*{}_*$-квазидетерминант*

*то мы снова получим истинное утверждение.*

**Теорема 3.9.**

$$(mA)^{-1_*{}^*}=A^{-1_*{}^*}m^{-1}\tag{3.18}$$

$$(Am)^{-1_*{}^*}=m^{-1}A^{-1_*{}^*}\tag{3.19}$$

[7]Мы можем дать подобное доказательство для $\binom{a}{b}$**-${}^*{}_*$-квазидетерминанта**. Однако мы можем записать соответствующие утверждения, опираясь на принцип двойственности. Так, если прочесть равенство (3.14) справа налево, то мы получим равенство

$$^j\det\left(A,{}^*{}_*\right)_i={}^jA_i-{}^{[j]}A_i{}^*{}_*\left({}^{[j]}A_{[i]}\right)^{-1^*{}_*}{}^*{}_*{}^jA_{[i]}$$

$$^j\det\left(A,{}^*{}_*\right)_i={}^jA_i-{}^{[j]}A_i{}^*{}_*\mathcal{H}\det\left({}^{[j]}A_{[i]},{}^*{}_*\right){}^*{}_*{}^jA_{[i]}$$

*Доказательство.* Мы докажем равенство (3.18) индукцией по размеру матрицы.

Для $1 \times 1$ матрицы утверждение очевидно, так как

$$(mA)^{-1_*{}^*} = ((mA)^{-1}) = (A^{-1}m^{-1}) = (A^{-1})\,m^{-1} = A^{-1_*{}^*}m^{-1}$$

Допустим утверждение справедливо для $(n-1)\times(n-1)$ матрицы. Тогда из равенства (3.1) следует

$$\begin{aligned}({}^I((mA)^{-1_*{}^*})_J)^{-1_*{}^*} =& {}^J(mA)_I - {}^J(mA)_{[I]*}{}^*\left({}^{[J]}(mA)_{[I]}\right)^{-1_*{}^*}{}_*{}^{*[J]}(mA)_I\\ =& m^JA_I - m\ {}^JA_{[I]*}{}^*\left({}^{[J]}A_{[I]}\right)^{-1_*{}^*}\ m^{-1}{}_*{}^*m\ {}^{[J]}A_I\\ =& m\ {}^JA_I - m\ {}^JA_{[I]*}{}^*\left({}^{[J]}A_{[I]}\right)^{-1_*{}^*}{}_*{}^{*[J]}A_I\end{aligned}$$

$$({}^I((mA)^{-1_*{}^*})_J)^{-1_*{}^*} = m\ {}^I(A^{-1_*{}^*})_J \tag{3.20}$$

Из равенства (3.20) следует равенство (3.18). Аналогично доказывается равенство (3.19). □

**Теорема 3.10.** *Пусть*

$$A = \begin{pmatrix}1 & 0\\ 0 & 1\end{pmatrix} \tag{3.21}$$

*Тогда*

$$A^{-1_*{}^*} = \begin{pmatrix}1 & 0\\ 0 & 1\end{pmatrix} \tag{3.22}$$

$$A^{-1^*{}_*} = \begin{pmatrix}1 & 0\\ 0 & 1\end{pmatrix} \tag{3.23}$$

*Доказательство.* Из (3.8) и (3.11) очевидно, что ${}^1(A^{-1_*{}^*})_1 = 1$ и ${}^2(A^{-1_*{}^*})_2 = 1$. Тем не менее выражение для ${}^2(A^{-1_*{}^*})_1$ и ${}^1(A^{-1_*{}^*})_2$ не может быть определено из (3.9) и (3.10) так как ${}^2A_1 = {}^1A_2 = 0$. Мы можем преобразовать эти выражения. Например

$$\begin{aligned}{}^2(A^{-1_*{}^*})_1 &= ({}^1A_2 - {}^1A_1({}^2A_1)^{-1}\ {}^2A_2)^{-1}\\ &= ({}^1A_1(({}^1A_1)^{-1}\ {}^1A_2 - ({}^2A_1)^{-1}\ {}^2A_2))^{-1}\\ &= (({}^2A_1)^{-1}\ {}^1A_1({}^2A_1({}^1A_1)^{-1}\ {}^1A_2 - {}^2A_2))^{-1}\\ &= ({}^1A_1({}^2A_1({}^1A_1)^{-1}\ {}^1A_2 - {}^2A_2))^{-1}\ {}^2A_1\end{aligned}$$

Мы непосредственно видим, что ${}^2(A^{-1_*{}^*})_1 = 0$. Таким же образом мы можем найти, что ${}^1(A^{-1_*{}^*})_2 = 0$. Это завершает доказательство (3.22).

Равенство (3.23) следует из (3.22), теоремы 3.7 и симметрии матрицы (3.21). □

## 4. Приводимое бикольцо

Если мы рассматриваем бикольцо матриц над полем $F$, то из коммутативности произведения в поле следует

$$A_*{}^*B = (A_a^cB_c^b) = (B_c^bA_a^c) = B^*{}_*A \tag{4.1}$$

**Определение 4.1.** **Приводимое бикольцо** - это бикольцо, в котором выполняется **условие приводимости произведений** (4.1). □

**Теорема 4.2.**

$$(A_*{}^*B)^T = B^T{}_*{}^*A^T \tag{4.2}$$

$$(A^*{}_*B)^T = B^T{}^*{}_*A^T \tag{4.3}$$

*в приводимом бикольце.*

*Доказательство.* Из (2.6) и (4.1) следует, что

$$(A_*{}^*B)^T = A^T{}^*{}_*B^T = B^T{}_*{}^*A^T$$

Мы можем доказать (4.3) аналогичным образом. □

**Теорема 4.3.**

$$A^{-1^*{}_*} = A^{-1_*{}^*} \tag{4.4}$$

*в приводимом бикольце.*

*Доказательство.* Из (2.22) и (4.1) следует, что

$$\delta = A_*{}^*A^{-1_*{}^*} = A^{-1_*{}^*}{}^*{}_*A \tag{4.5}$$

(4.4) следует из сравнения (4.5) и (2.23). □

## 5. Квазидетерминант над полем

**Определение 5.1.** Для матриц, определённых над полем, мы определим функцию, называемую **определитель матрицы**,

$$\det() = 1 \tag{5.1}$$

$$\det a = \sum_a (-1)^{|a|+|b|} a^b_a \det a^{[a]}_{[b]} \tag{5.2}$$

□

**Теорема 5.2.**

$$\det (a,{}_*{}^*)^a_b = (-1)^{|a|+|b|} \frac{\det a}{\det a^{[b]}_{[a]}} \tag{5.3}$$

*для матриц, определённых над полем.*

*Доказательство.* Мы проведём доказательство индукцией по $n$. Так как умножение в поле коммутативно, мы можем уточнить выражение (3.15)

$$\det (a,{}_*{}^*)^a_b = a^a_b - \frac{a^c_b a^a_e}{\det \left(a^{[b]}_{[a]},{}_*{}^*\right)^c_e} \tag{5.4}$$

При $n = 1$ мы можем доказать утверждение теоремы непосредственной проверкой.

Пусть утверждение теоремы справедливо при $n = k - 1$. Тогда это утверждение справедливо для всех миноров $k \times k$ матрицы $a$. Мы подставим (5.3) в равенство (5.4)

$$\begin{aligned}
\det (a, {}_*{}^*)^a_b = a^a_b &- \sum_{|c|<|a|} \sum_{|e|<|b|} \frac{a^c_b a^a_e}{(-1)^{|c|+|e|} \frac{det\ a[^b_a]}{det\ a[^{b,e}_{a,c}]}} \\
&+ \sum_{|c|>|a|} \sum_{|e|<|b|} \frac{a^c_b a^a_e}{(-1)^{|c|+|e|} \frac{det\ a[^b_a]}{det\ a[^{b,e}_{a,c}]}} \\
&+ \sum_{|c|<|a|} \sum_{|e|>|b|} \frac{a^c_b a^a_e}{(-1)^{|c|+|e|} \frac{det\ a[^b_a]}{det\ a[^{b,e}_{a,c}]}} \\
&- \sum_{|c|>|a|} \sum_{|e|>|b|} \frac{a^c_b a^a_e}{(-1)^{|c|+|e|} \frac{det\ a[^b_a]}{det\ a[^{b,e}_{a,c}]}}
\end{aligned}$$

$$\begin{aligned}
\det (a, {}_*{}^*)^a_b = a^a_b \quad &- \sum_{|c|<|a|} \sum_{|e|<|b|} \frac{(-1)^{|a|+|b|+|e|+|c|} a^c_b a^a_e det\ a[^{b,e}_{a,c}]}{(-1)^{|a|+|b|} det\ a[^b_a]} \\
&+ \sum_{|c|>|a|} \sum_{|e|<|b|} \frac{(-1)^{|a|+|b|+|e|+|c|} a^c_b a^a_e det\ a[^{b,e}_{a,c}]}{(-1)^{|a|+|b|} det\ a[^b_a]} \\
&+ \sum_{|c|<|a|} \sum_{|e|>|b|} \frac{(-1)^{|a|+|b|+|e|+|c|} a^c_b a^a_e det\ a[^{b,e}_{a,c}]}{(-1)^{|a|+|b|} det\ a[^b_a]} \\
&- \sum_{|c|>|a|} \sum_{|e|>|b|} \frac{(-1)^{|a|+|b|+|e|+|c|} a^c_b a^a_e det\ a[^{b,e}_{a,c}]}{(-1)^{|a|+|b|} det\ a[^b_a]}
\end{aligned} \tag{5.5}$$

Подставив равенство

$$\det a[^b_c] = \begin{cases} \left( \sum\limits_{|e|<|b|} - \sum\limits_{|e|>|b|} \right) (-1)^{|a|+|e|} a^a_e det\ a[^{b,e}_{a,c}] & c > a \\ \left( \sum\limits_{|e|>|b|} - \sum\limits_{|e|<|b|} \right) (-1)^{|a|+|e|} a^a_e det\ a[^{b,e}_{a,c}] & c < a \end{cases}$$

в равенство (5.5), мы получим

$$\det (a, {}_*{}^*)^a_b = a^a_b + \left( \sum_{|c|<|a|} + \sum_{|c|>|a|} \right) \frac{(-1)^{|b|+|c|} a^c_b det\ a[^b_c]}{(-1)^{|a|+|b|} det\ a[^b_a]}$$

Сложив дроби, мы получим

$$\det (a, {}_*{}^*)^a_b = (-1)^{|a|+|b|} \sum_{c \in M} \frac{(-1)^{|b|+|c|} a^c_b \det a[^b_c]}{det\ a[^b_a]} \tag{5.6}$$

Подставив (5.2) в (5.6), мы докажем, что утверждение теоремы справедливо при $n = k$. □

## 6. Список литературы

## 7. Предметный указатель

# 8. Специальные символы и обозначения